\documentclass[12pt, twoside, leqno]{article}
\usepackage{amsmath,amsthm}
\usepackage{amssymb,latexsym}
\usepackage{enumerate}

\newtheorem{theorem}{Theorem}[section]
\newtheorem{corollary}[theorem]{Corollary}
\newtheorem{lemma}[theorem]{Lemma}
\newtheorem{proposition}[theorem]{Proposition}

\theoremstyle{definition}

\numberwithin{equation}{section}

\def\Z{\mathbb{Z}}
\def\F{\mathbb{F}}
\def\R{\mathbb{R}}

\def\eps{\varepsilon}
\def\qed{$\blacksquare$}

\makeatletter
\def\pmod#1{\allowbreak\mkern10mu({\operator@font mod}\,\,#1)}
\def\pod#1{\allowbreak\mkern10mu(#1)}
\makeatother

\begin{document}

\baselineskip=17pt


\title{On the constant in Burgess' bound for the number of consecutive residues or non-residues}

\author{Kevin J. McGown\\
Department of Mathematics\\
University of California, San Diego\\
9500 Gilman Drive\\
La Jolla, CA 92093-0112\\
E-mail:  kmcgown@math.ucsd.edu
}

\date{}

\maketitle


\renewcommand{\thefootnote}{}



\renewcommand{\thefootnote}{\arabic{footnote}}
\setcounter{footnote}{0}


\begin{abstract}
We give an explicit version of a result due to D. Burgess.
Let $\chi$ be a non-principal Dirichlet character modulo a prime $p$.
We show that the maximum number of consecutive integers
for which $\chi$ takes on a particular value is less than
$\left\{\frac{\pi e\sqrt{6}}{3}+o(1)\right\}p^{1/4}\log p$,
where the $o(1)$ term is given explicitly.
\end{abstract}

%
%
%
%
%
%
%
%
%


\section{Introduction}\label{S:1}

Let $\chi$ be a non-principal Dirichlet character to the prime modulus $p$.
In 1963, D.~Burgess showed (see~\cite{bib:burgess.1963}) that the maximum number of consecutive integers for which $\chi$ takes
on any particular value is $O(p^{1/4}\log p)$.
This still constitutes the best known asymptotic upper bound on this quantity.
However, in some applications, one needs a more explicit result.
Following the general lines of his original argument and
making careful estimates throughout,
we prove an explicit version of Burgess' theorem (see Theorem~\ref{T:3} and Corollary~\ref{C:1}),
thereby obtaining the following:
\begin{theorem}\label{T:1}
If $\chi$ is any non-principal Dirichlet character to the prime modulus $p$
 which is constant on $(N,N+H]$, then
$$
  H<\left\{\frac{\pi e\sqrt{6}}{3}+o(1)\right\}
  p^{1/4}\log p
  \,.
$$
\end{theorem}
We note that the constant $(\pi e\sqrt{6})/3$ is approximately $6.97$.
As we have an explicit bound on the $o(1)$ term when $p$ is large,
we are able to obtain the following result which is more useful in applications:
\begin{theorem}\label{T:2}
If $\chi$ is any non-principal Dirichlet character to the prime modulus $p$
which is constant on $(N,N+H]$, then
$$
  H
  <
  \begin{cases}
  7.06\,p^{1/4}\log p
  \,,\;\;
  &
  \text{for }p\geq 5\cdot 10^{18}
  \\[0.5ex]
  7\,p^{1/4}\log p
  \,,\;\;
  &
  \text{for }p\geq 5\cdot10^{55}
  \end{cases}
  \,.  
$$
\end{theorem}

For the special case of $N=0$, which amounts to giving a bound on the smallest non-residue of $\chi$
(i.e., the smallest $n$ such that $\chi(n)\neq 1$),
K.~Norton proves a result analogous to Theorem~\ref{T:2}
which holds for all $p$ with a constant of $4.7$ (see~\cite{bib:norton.1971}).
In addition, a result for arbitrary $N$, similar to the one given in Theorem~\ref{T:2} is stated,
but not proved in~\cite{bib:norton.1973}.
R.~Hudson (see~\cite{bib:hudson.1974}) cites a result slightly improving the
one stated in~\cite{bib:norton.1973} to appear in a future paper,
but the present author cannot locate the purported proof.
It seems a worthwhile endeavor to put down such a proof as it is possible that some
authors avoid using the result in~\cite{bib:norton.1973} due to the lack of proof
(see, for example~\cite{bib:hummel.2003}), while others (see~\cite{bib:hudson.1974})
use the result for further derivations.
To our knowledge, this is the first proof to appear in the literature which makes
the constant in Burgess' theorem explicit.
\\


It is perhaps useful here to comment briefly on the connection between Dirichlet characters and power residues.
Fix an integer $k\geq 2$.
We say that $n\in\Z$ is a $k$-th power residue modulo $p$ if 
$(n,p)=1$ and the equation $x^k\equiv n\pmod{p}$ is soluble in $x$.
Suppose $\chi$ is any Dirichlet character modulo $p$ of order $(k,p-1)$.
One can easily show that
$\chi(n)=1$ if and only if $n$ is a $k$-th power residue modulo $p$.
Here we might as well assume $(k,p-1)>1$, or else every integer is a $k$-th power residue modulo $p$
and the only such $\chi$ is the principal character.
If we denote by $C_p=(\Z/p\Z)^\star$ the multiplicative group consisting
of the integers modulo $p$ and by $C_p^k$ the subgroup of $k$-th powers modulo $p$, then the value of
$\chi(n)$ determines to which coset of $C_p/C_p^k$ the integer $n$ belongs.
In light of this, theorems~\ref{T:1} and \ref{T:2} also give estimates
(which are the best known)
on the maximum number of consecutive integers that belong to
a given coset of $C_p/C_p^k$.
\\

We should also mention that Burgess' well-known character sum estimate (see~\cite{bib:burgess.1962}) gives a bound
on the quantity (in the title of the paper) of $O(p^{1/4+\eps})$.
However, the constant associated to the $O$-symbol depends on $\eps$ and hence,
although there are explicit versions of Burgess' character sum estimate available (see~\cite{bib:IW}),
theorems~\ref{T:1} and~\ref{T:2} would not follow from
this.
\\

The main idea behind Burgess' proof is to combine upper and lower bounds for the sum:
$$
 S(\chi,h,r):= \sum_{x=0}^{p-1}
  \left|
  \sum_{m=1}^h
  \chi(x+m)
  \right|^{2r}
$$
In Lemma~\ref{L:1C} of \S\ref{S:2} we give an upper bound for $S(\chi,h,r)$ in terms of $r$ and $h$.
In Proposition~\ref{P:1} of \S\ref{S:3} we give a lower
bound on $S(\chi,h,r)$ in terms of $h$ and $H$, under some additional hypotheses on $H$.
Combining these results, we obtain an upper bound on $H$ in terms of $r$ and $h$ under the same hypotheses;
this result is also given as part of Proposition~\ref{P:1}.
Then, in \S\ref{S:4} we prove our main result (see Theorem~\ref{T:3}) by invoking Proposition~\ref{P:1} with a careful choice of parameters.
Finally, by performing some simple numerical computations, we show that that the extra hypothesis on $H$ can be
dropped when $p$ is large enough (see Corollary~\ref{C:1}); theorems~\ref{T:1} and \ref{T:2} will then follow immediately.


\section{An Upper Bound on $S(\chi,h,r)$}\label{S:2}

The following character sum estimate was first given by A. Weil, as a consequence of his
deep work on the Riemann hypothesis for function fields (see~\cite{bib:weil}).
It is also proved as Theorem~2C' in \cite{bib:schmidt} 
using an elementary method due to S.~Stepanov (see~\cite{bib:stepanov}),
which was later extended by both E.~Bombieri (see~\cite{bib:bombieri.1973}) and W. Schmidt (see~\cite{bib:schmidt.1973}).
\begin{lemma}\label{L:Weil}
Let $\chi$ be a non-principal Dirichlet character to the prime modulus $p$,
having order $n$.
Let $f(x)\in\Z[x]$ be a polynomial with $m$ distinct roots
which is not an $n$-th power in $\F_p[x]$,
where $\F_p$ denotes the finite field with $p$ elements.
Then
$$
  \left|
  \sum_{x\in\F_p}
  \chi(f(x))
  \right|
  \leq
  (m-1)\;
  p^{1/2}
  \,.
$$
\end{lemma}
The next lemma is a slight improvement over
Lemma 2 in~\cite{bib:burgess.1962} which gives
an upper bound on $S(\chi,h,r)$.  The proof is not difficult if we grant ourselves
Lemma~\ref{L:Weil}.
\begin{lemma}\label{L:1C}
Suppose $\chi$ is any non-principal Dirichlet character to the prime modulus $p$.
If $r,h\in\Z^+$, then
$$
 S(\chi,h,r)
  <
  \frac{1}{4}(4r)^rp h^r
  +
  (2r-1)p^{1/2}h^{2r}
  \,.
$$  
\end{lemma}

\noindent\textbf{Proof.}
First we claim that we may assume, without loss of generality,
that \mbox{$r<h<p$}.  We commence by observing
that $h=p$ implies $S(\chi,h,r)=0$, in which case there is nothing to prove.
We see that $h>p$ implies $S(\chi,h-p,r)=S(\chi,h,r)$, which allows us
to inductively bring $h$ into the range $0<h<p$.
Additionally, we notice that if $h\leq r$, then the theorem
is trivial since in this case we would have
$S(\chi,h,r)\leq h^{2r}p\leq (hr)^rp$.
This establishes the claim.

Now, to begin the proof proper, 
we observe that
$$
  S(\chi,h,r)
  =
  \sum_{1\leq m_1,\dots,m_{2r}\leq h}\;\;
  \sum_{x=0}^{p-1}
  \chi(x+m_1)\dots\chi(x+m_r)\overline{\chi}(x+m_{r+1})\dots\overline{\chi}(x+m_{2r})
  \,.
$$
Define
$$
  \mathcal{M}:=\{\mathbf{m}=(m_1,\dots,m_{2r})\mid 1\leq m_1,\dots,m_{2r}\leq h\}
  \,.
$$
We can rewrite the above as
$$
  S(\chi,h,r)=
  \sum_{\mathbf{m}\in\mathcal{M}}
  \sum_{x\in\F_p}
  \chi(f_\mathbf{m}(x))
  \,,
$$
where
$$
  f_\mathbf{m}(x)
  =
  (x+m_1)\dots(x+m_r)(x+m_{r+1})^{n-1}(x+m_{2r})^{n-1}
  \,,
$$
and $n$ denotes the order of $\chi$.
%
If $f_\mathbf{m}(x)$ is not an $n$-th power mod $p$, then by Lemma~\ref{L:Weil}
we have
$$
  \left|
  \sum_{x\in\F_p}
  \chi(f_\mathbf{m}(x))
  \right|
  \leq
  (2r-1)\sqrt{p}
  \,.
$$
Otherwise, we must settle for the trivial bound of $p$.

It remains to count the number of exceptions -- that is, the number of $\mathbf{m}\in\mathcal{M}$ such that
$f_\mathbf{m}(x)$ is an $n$-th power mod $p$.
A little care is required here -- as an example, if $r=n=3$ and $p\geq 5$, then the vectors
$\mathbf{m}=(1,2,3,1,2,3)$ and $\mathbf{m}=(1,1,1,2,2,2)$ are both exceptions, but the way in which they
arise is slightly different;
as $r$ gets larger compared to $n$, the situation only gets worse.
In light of this difficulty, we will actually count (as Burgess does in~\cite{bib:burgess.1963}) the number of
$\mathbf{m}=(m_1,\dots,m_{2r})\in\mathcal{M}$ such that each $m_j$ is repeated at least
once.

We let $u$ denote the number of distinct $m_j$ (so that $u\leq r<h$)
and denote by \mbox{$1=j_1<j_2<\dots<j_u\leq 2r$} the indices corresponding to the first
occurrence of each of the $u$ values among the $m_j$.
The number of ways to choose the $j_k$ is bounded by $\binom{2r-1}{u-1}$,
and there are at most $h$ choices for each $m_{j_k}$ while the remaining $m_{j}$ are restricted
to at most $u$ values.
In light of all this, we find that 
the number of exceptions 
is bounded above by
\begin{eqnarray*}
  \sum_{u=1}^r\binom{2r-1}{u-1}h^uu^{2r-u}
  &\leq&
  (hr)^r\sum_{u=1}^r\binom{2r-1}{u-1}\left(\frac{u}{h}\right)^{r-u}
 \;\leq\;
  (hr)^r\sum_{u=1}^r\binom{2r-1}{u-1}
  \,.
\end{eqnarray*}
Finally, to complete the proof, we observe
$$
  (hr)^r\sum_{u=1}^r\binom{2r-1}{u-1}
   \;=\;
  (hr)^r2^{2r-2}
  \;=\;
  \frac{1}{4}
  (4rh)^r
  \,.
  \;\;
  \text{\qed}
$$


\section{A Lower Bound on $S(\chi,h,r)$}\label{S:3}

In obtaining the desired lower bound, the idea is to locate a large number of intervals on which $\chi$ is constant.
The next two lemmas will be useful in accomplishing this end.
The following lemma makes the error term in Lemma 3 of~\cite{bib:burgess.1963} explicit
and improves the main constant from
$1-\pi^2/12\approx 0.178$ to $3/\pi^2\approx 0.304$.
\begin{lemma}\label{L:3B}
  Let $X\geq 7$.  If $a,b\in\Z$ are coprime with $a\geq 1$, then there are at least
  $$
    X^2
   \left(\frac{3}{\pi^2}-\frac{\log X}{2X}-\frac{1}{X}-\frac{1}{2X^2}\right)
  $$
  distinct numbers of the form
  $$
    \frac{at+b}{q}
  $$
  where $0\leq t<q\leq X$.
\end{lemma}


\noindent\textbf{Proof.}
As in \cite{bib:burgess.1963}, we observe that
$
  \#\{q^{-1}(at+b) \mid 0\leq t<q\leq X\}
$
is bounded below by
\begin{eqnarray*}
  \sum_{q\leq X}\sum_{\substack{0\leq t<q\\(at+b,q)=1}}
  1
  &=&
  \sum_{q\leq X}
  \;
  \sum_{0\leq t<q}
  \;
  \sum_{m|(at+b,q)}
  \mu(m)
\end{eqnarray*}
Writing $q=rm$ allows us to rewrite the sum above as
\begin{eqnarray*}
  \sum_{m\leq X}\mu(m)
  \sum_{r\leq X/m}
  \sum_{\substack{0\leq t<rm\\at\equiv -b\pod m}}
  1
  \,.
\end{eqnarray*}

Since $(a,b)=1$, the congruence $at\equiv -b\pmod m$ has a solution if and only if $(m,a)=1$.
Therefore we can rewrite our sum in the following way:
\begin{eqnarray*}
  \sum_{\substack{m\leq X\\ (m,a)=1}}
  \mu(m)
  \sum_{r\leq X/m}
  \sum_{\substack{0\leq t<rm\\at\equiv -b\pod m}}
  1
  &=&
  \sum_{\substack{m\leq X\\ (m,a)=1}}
  \mu(m)
  \sum_{r\leq X/m}
  r
\end{eqnarray*}
A careful lower estimate of the sum on the right-hand side above
will give the desired result.
Using the identity
$$
  \sum_{r\leq Y}r=\frac{Y^2}{2}+\frac{Y}{2}\theta_Y
  \,,
  \quad
  \theta_Y\in[-1,1]
  \,,
$$
which holds for $Y>0$, we obtain
\begin{eqnarray}\label{E:BigMu}
  &&
  \sum_{\substack{m\leq X\\ (m,a)=1}}
  \mu(m)
  \sum_{r\leq X/m}
  r
  \\
  \nonumber
  &&\qquad\qquad\qquad
  =
  \;\;
  \frac{X^2}{2}
  \sum_{\substack{m\leq X\\ (m,a)=1}}
  \frac{\mu(m)}{m^2}
  \;+\;
  \frac{X}{2}
    \sum_{\substack{m\leq X\\ (m,a)=1}}
    \frac{\mu(m)}{m}\theta_{X/m}
    \,.  
\end{eqnarray}

Let $\zeta(s)$ denote the Riemann zeta function.
When $s>1$, we have
\begin{eqnarray*}
\sum_{\substack{m=1\\(m,a)=1}}^\infty
\mu(m) m^{-s}
\;=\;
\zeta(s)^{-1}\prod_{p | a}(1-p^{-s})^{-1}
\;\geq\;
\zeta(s)^{-1}
\,,
\end{eqnarray*}
and the tail of the series is bounded in absolute value by
\begin{eqnarray*}
\sum_{m>X}m^{-s}
\leq
\frac{1}{X^{s}}+\frac{1}{(s-1)}\cdot \frac{1}{X^{s-1}}
\,;
\end{eqnarray*}
therefore
$$
  \sum_{\substack{m\leq X\\(m,a)=1}}
  \mu(m) m^{-s}
  \geq
  \zeta(s)^{-1}
  -\frac{1}{X^{s}}-\frac{1}{(s-1)}\cdot \frac{1}{X^{s-1}}
  \,.
$$
Setting $s=2$ gives
$$
  \sum_{\substack{m\leq X\\(m,a)=1}}
  \frac{\mu(m)}{m^2}
  \geq
  \zeta(2)^{-1}
  -\frac{1}{X^2}-\frac{1}{X}
  \,.
$$

Now we deal with the second sum on the right-hand side of (\ref{E:BigMu});
we have
$$
\vrule\;
  \sum_{\substack{m\leq X\\(m,a)=1}}^{\phantom{N}}
  \frac{\mu(m)}{m}\theta_{X/m}
\;\vrule
  \leq
  \sum_{m\leq X}\frac{1}{m}
  \leq
  1+\log X
  \,.
$$
Summarizing, we have shown
\begin{eqnarray*}
  \vrule\;
  \sum_{\substack{m\leq X\\ (m,a)=1}}^{\phantom{N}}
  \mu(m)
  \sum_{r\leq X/m}
  r
  \;\vrule
  &\geq&
  \frac{X^2}{2}
  \left(\frac{1}{\zeta(2)}-\frac{1}{X^2}-\frac{1}{X}\right)
  -
  \frac{X}{2}
  \left(1+\log X\right)
  \\
  &=&
  X^2
  \left(\frac{1}{2\zeta(2)}-\frac{\log X}{2X}-\frac{1}{X}-\frac{1}{2X^2}\right)
  \,.
\end{eqnarray*}
In light of the fact that $\zeta(2)=\pi^2/6$, we have arrived at the desired conclusion.
The reader may worry why we failed to use the hypothesis that $X\geq 7$.
This hypothesis is not necessary for the truth of the conclusion,
but we include it nonetheless to ensure that our estimate gives a positive number.
\qed

\vspace{1ex}
Finally we will require Dirichlet's Theorem in Diophantine approximation; see, for example,
Theorem 1 in Chapter 1 of~\cite{bib:cassels}.
\begin{lemma}\label{L:2}
  Let $\theta,\,A\in\R$ with $A>1$.
  Then there exists $a,b\in\Z$ with $(a,b)=1$ such that
  $$
    0<a<A
    \,,\quad
    |a\theta-b|\leq A^{-1}
    \,.
  $$  
\end{lemma}

We are now ready to give our lower bound on $S(\chi,h,r)$.


\begin{proposition}\label{P:1}
  Let $h,r\in\Z^+$.
  Suppose $\chi$ is a non-principal Dirichlet character to the prime modulus $p$ which
  is constant on $(N,N+H]$ and such that 
  $$
    14h\leq H\leq (2h-1)^{1/3}p^{1/3}
    \,.
  $$  
  If we set $X:=H/(2h)\geq 7$, then
  $$
    S(\chi,h,r)\geq \left(\frac{3}{\pi^2}\right)X^2 h^{2r+1} f(X)
    \,,
  $$
  where
  $$
    f(X)=1-\frac{\pi^2}{3}\left(\frac{\log X}{2X}+\frac{1}{X}+\frac{1}{2X^2}\right)
    \,,
  $$
  and therefore
  \begin{equation}
  \nonumber
  H
  <
  \frac{2\pi h}{\sqrt{3f(X)}}
  \,
  p^{1/4}
  \left[
  \frac{1}{4h}
  \left(\frac{4r}{h}\right)^{r}p^{1/2}
  +
  \left(\frac{2r-1}{h}\right)
  \right]^{1/2}
    \,.
  \end{equation}
  Note:  $f(X)$ is positive and increasing on $[7,\infty)$ and $f(X)\to 1$ as $X\to\infty$.
%
\end{proposition}

\noindent\textbf{Proof.}
Following the argument given in \cite{bib:burgess.1963},
we define the real interval
$$
  I(q,t):=
  \left(\frac{N+pt}{q},\;\; \frac{N+H+pt}{q}\right]
  \,,
$$
for $0\leq t<q\leq X$.
We take note of two important properities of $I(q,t)$, which we will use later.
First, the length of $I(q,t)$ is $H/q\geq H/X=2h$.  Second, $\chi$ is constant on $I(q,t)$;
this is because for any $z\in I(q,t)$ we have $\chi(qz-pt)=\zeta$
and hence $\chi(z)=\overline{\chi}(q)\zeta$.
We are interested in locating a large number of non-overlapping intervals of this form.

By Lemma~\ref{L:2},
there exists coprime $a,b\in\Z$ such that
$1\leq a\leq H$ and
\begin{equation}\label{E:I1}
  |aNp^{-1}-b|\leq 1/H
  \,.
\end{equation}
One shows
that if $I(q_1,t_1)$ and $I(q_2,t_2)$ overlap, then
\begin{equation}\label{E:I2}
  |Np^{-1}(q_1-q_2)+t_2 q_1- t_1 q_2|
  <
  p^{-1}XH
  \,.
\end{equation}
Equations
(\ref{E:I1}) and (\ref{E:I2}) yield
$$
  \left|
  \frac{b}{a}
  (q_1-q_2)+t_2 q_1-t_1 q_2
  \right|
  <
  \frac{XH}{p}+\frac{|q_1-q_2|}{Ha}
  \leq
  \frac{XH}{p}+\frac{X}{Ha}
  =
  \frac{H^2 a+p}{2ahp}
  \,.
$$
But since $a\leq H$ and $H^3\leq (2h-1)p$ by hypothesis, we have
$$
  \frac{H^2 a+p}{2ahp}
  \leq
  \frac{H^3+p}{2ahp}
  \leq
  \frac{1}{a}
  \,.
$$
Hence
$$
  \left|
  \frac{b}{a}
  (q_1-q_2)+t_2 q_1-t_1 q_2
  \right|
  <
  \frac{1}{a}
  \,,
$$
and it follows that
$I(q_1,t_1)$ and $I(q_2,t_2)$ can only overlap if
$$
  \frac{a t_1+b}{q_1}=\frac{a t_2+b}{q_2}
  \,.
$$
Invoking Lemma \ref{L:3B}, we find that there will be at least $(3/\pi^2)X^2 f(X)$ disjoint intervals $I(q,t)$ of the given form.
\\

Having located the desired intervals, we are ready to give a lower estimate
for $S(\chi,h,r)$.
Let $z(q,t)$ denote the smallest integer in $I(q,t)$.
Since $I(q,t)$ has length at least $2h$,
the integers $z(q,t)+n+m$, for $n=0,\dots,h-1$ and $m=1,\dots,h$ are distinct elements of $I(q,t)$.
Moreover, as $q,t$ run through the values selected by Lemma~\ref{L:3B}, the $I(q,t)$ are disjoint.
Now, using the fact that $\chi$ is constant on each $I(q,t)$, one obtains the following bound for $S(\chi,h,r)$:
\begin{eqnarray*}
  \sum_{x=0}^{p-1}
  \left|
  \sum_{m=1}^h\chi(x+m)\right|^{2r}
  &\geq&
  \sum_{q,t}\sum_{n=0}^{h-1}
  \left|
  \sum_{m=1}^{h}\chi(z(q,t)+n+m)
  \right|^{2r}
  \\
  &=&
  \sum_{q,t}  
  \sum_{n=0}^{h-1}h^{2r}
  \\
  &=&
  h^{2r+1}\sum_{q,t} 1
  \\
  &\geq&
  \left(\frac{3}{\pi^2}\right)X^2 h^{2r+1} f(X)
  \,.
\end{eqnarray*}

Now we combine this lower bound on
$S(\chi,h,r)$ with the upper bound given in Lemma~\ref{L:1C}
to obtain
$$
  \left(\frac{3}{\pi^2}\right)
  \left(\frac{H}{2h}\right)^2 h^{2r+1} f(X)
  <
  \frac{1}{4}(4r)^{r}ph^r + (2r-1)p^{1/2}h^{2r}
  \,,
$$
which implies
\begin{eqnarray*}
  H^2
  &<&
  \frac{4\pi^2h^2}{3f(X)}
  \,
  p^{1/2}
  \left[
  \frac{1}{4h}
  \left(\frac{4r}{h}\right)^{r}p^{1/2}
  +
  \left(\frac{2r-1}{h}\right)
  \right]
  \,.
\end{eqnarray*}
(We have used the fact that $f(X)>0$ for $X\geq 7$ in order to divide both sides
by $f(X)$ and preserve the inequality.)
Taking the square root of both sides yields the result.
\qed


\section{The Main Result}\label{S:4}

\begin{theorem}\label{T:3}
  Suppose $\chi$ is a non-principal Dirichlet character to the prime modulus $p\geq 5\cdot 10^4$
  which is constant on $(N,N+H]$.
  If
  $H\leq (2e^2\log p-3)^{1/3}p^{1/3}$,
  then
  $$
  H<
  C\,g(p)\cdot p^{1/4}\log p\cdot
  $$
  where
  $$
    C=\frac{\pi e\sqrt{6}}{3}\approx 6.97266
  $$
  and
  $g(p)\to 1$ as $p\to\infty$.  In fact,
  $$
    g(p)=
    \sqrt{
    f\left(\frac{Cp^{1/4}}{2e^2}\right)^{-1}
    \left(
    1+\frac{1}{\log p}
    \right)
    }
    \,,
  $$
  where $f(X)$ is defined in Proposition~\ref{P:1}.
  Note that $g(p)$ is positive and decreasing for $p\geq 5\cdot 10^4$.
\end{theorem}
Before launching the proof of Theorem~\ref{T:3}, we will establish the following:
\begin{lemma}\label{L:logs}
  Let $p\geq 3$ be an integer.
  Suppose that $A,B>0$ are real numbers such that $h=\lfloor A\log p\rfloor$ and $r=\lfloor B\log p\rfloor$ are positive integers with $2r+1\leq h$.
  Then
  $$
    A\geq 4B\cdot\exp\left(\frac{1}{2B}\right)
    \;\;\Longrightarrow\;\;
    \frac{1}{2h}\left(\frac{4r}{h}\right)^r\leq \frac{1}{Ap^{1/2}\log p}
    \,.
  $$

\end{lemma}

\noindent\textbf{Proof.}
By convexity, $\log t\geq (2\log 2)(t-1)$ for all $t\in[1/2,1]$ and thus
$$
  \log\left(\frac{h}{h+1}\right)\geq
  \frac{-2\log 2}{h+1}
  \geq
  \frac{-\log 2}{r+1}
  \,.
$$
This implies
$$
  \frac{1}{2}\leq\left(\frac{h}{h+1}\right)^{r+1}
$$
and therefore
$$
  \frac{1}{2h}\left(\frac{4r}{h}\right)^{r}
  \leq
  \frac{1}{h+1}
  \left(\frac{4r}{h+1}\right)^r
  \leq
  \frac{1}{A\log p}
  \left(\frac{4B}{A}\right)^r
  \,.
$$
Hence to obtain the desired implication, is suffices to show
\begin{equation}
\nonumber
  \left(\frac{4B}{A}\right)^r\leq p^{-1/2}
  \,.
\end{equation}
Taking logarithms, this is equivalent to
$$
  r\log \left(\frac{4B}{A}\right)\leq -\frac{1}{2}\log p
  \,,
$$
which follows from inequality
$$
  B\log \left(\frac{4B}{A}\right)\leq -\frac{1}{2}
  \,,
$$
which is true by hypothesis.
\qed

\vspace{2ex}

\noindent\textbf{Proof of Theorem~\ref{T:3}.}
We will suppose $H\geq Cp^{1/4}\log p$, or else there is nothing to prove.
Set $h=\lfloor A\log p\rfloor$ and $r=\lfloor B\log p\rfloor$,
where $A:=e^2$ and $B:=1/4$.  The constants $A$ and $B$ were chosen as to minimize the quantity
$AB$ subject to the constraint $A\geq 4B\exp\left(\frac{1}{2B}\right)$.
One easily checks that $14h\leq Cp^{1/4}\log p$ for our choices of $h$ and $C$, provided $p\geq 5\cdot 10^4$ and hence $14h\leq H$.
Also, we note that $H\leq(2h-1)^{1/3}p^{1/3}$ by hypothesis.
We apply Proposition~\ref{P:1} and adopt all notation relevant to its statement.
This gives:
  \begin{equation}\label{E:bound}
  H
  <
  \frac{2\pi h}{\sqrt{3f(X)}}
  \,
  p^{1/4}
  \left[
  \frac{1}{4h}
  \left(\frac{4r}{h}\right)^{r}p^{1/2}
  +
  \left(\frac{2r-1}{h}\right)
  \right]^{1/2}
  \end{equation}
  
In order for the quantity inside the square brackets above to remain bounded as $p$ gets large,
and moreover be as small as possible, we would like
\begin{equation}\label{E:mycond}
  \nonumber
  \frac{1}{4h}\left(\frac{4r}{h}\right)^rp^{1/2}\to 0
  \,.
\end{equation}
As the constants $A$ and $B$ were chosen to satisfy the conditions of Lemma~\ref{L:logs}
(the condition above was precisely the motivation for the lemma),
we have
$$
\frac{1}{2h}
  \left(\frac{4r}{h}\right)^{r}
  \leq
  \frac{1}{Ap^{1/2}\log p}
  \,.
$$
To give a clean bound on the the quantity $(2r-1)/h$ we notice that $2r\leq h+1$ implies
$$
  \frac{2r-1}{h}
  \leq\frac{2r}{h+1}
  \leq
  \frac{2B}{A}
  \,.
$$

Thus inequality (\ref{E:bound}) becomes
\begin{eqnarray*}
   H
   &<&
   \frac{2\pi A}{\sqrt{3f(X)}}\,p^{1/4}\log p
   \left[
    \frac{1}{2A\log p }+\frac{2B}{A}
   \right]^{1/2}
   \\
   &=&
    p^{1/4}\log p
   \left[
     \frac{8\pi^2AB}{3f(X)}\left(
    1+\frac{1}{4B\log p}
    \right)
   \right]^{1/2}
   \,.
\end{eqnarray*}
Now it is plain that the asymptotic constant in the above expression is directly proportional to $\sqrt{AB}$,
which motivates our choices of $A$ and $B$.
Plugging in the values of $A$ and $B$, we obtain:
\begin{eqnarray*}
   H
   &<&
   p^{1/4}\log p
   \left[
      \frac{2\pi^2e^2}{3f(X)}
  \left(
   1+\frac{1}{\log p}
   \right)
    \right]^{1/2}
       \\
   &=&
   \frac{e\pi\sqrt{6}}{3}
   \,
   p^{1/4}\log p
   \left[
      \frac{1}{f(X)}
  \left(
   1+\frac{1}{\log p}
   \right)
    \right]^{1/2}
\end{eqnarray*}
Finally, we note that we have an a priori lower bound on $X$; namely
$$
  X=\frac{H}{2h}
  \geq
  \frac{Cp^{1/4}\log p}{2A\log p}
  =
  \frac{C\,p^{1/4}}{2e^2}
  \,.
$$  
In light of the fact that $f(X)$ is increasing, this gives
$$
  f(X)^{-1}\leq f\left(\frac{C\,p^{1/4}}{2e^2}\right)^{-1}
  \,,
$$
and the result follows.
\qed

\begin{corollary}\label{C:1}
  If $\chi$ is a non-principal Dirichlet character to the prime modulus $p\geq 5\cdot 10^{18}$
  which is constant on $(N,N+H]$,
  then
  $$
  H<
  C\,g(p)\cdot p^{1/4}\log p
  \,,
  $$
  where $C$ and $g(p)$ are as in Theorem~\ref{T:3}.
\end{corollary}

\noindent\textbf{Proof.}
In order to apply Theorem~\ref{T:3}, which will give the result, it suffices to show that
$H\leq (2e^2\log p-3)^{1/3}p^{1/3}$.  By way of contradiction,
suppose $H>(2e^2\log p-3)^{1/3}p^{1/3}$.  In this case
we set $H=\lfloor (2e^2\log p -3)^{1/3}p^{1/3}\rfloor$,
and note that $\chi$ is clearly still constant on $(N,N+H]$ for smaller $H$.
We invoke Theorem~\ref{T:3} to conclude that
$H<Cg(p)p^{1/4}\log p$ where $Cg(p)\leq Cg(5\cdot 10^{18})< 7.06$.
Using the fact that $p\geq 5\cdot 10^{18}$, we have
$$
H<7.06 p^{1/4}\log p<(2e^2\log p-3)^{1/3}p^{1/3}-1<H
\,,
$$
which is
a contradiction.
\qed

\vspace{1ex}
It remains to derive theorems~\ref{T:1} and \ref{T:2}.
Theorem~\ref{T:1} follows immediately from Corollary~\ref{C:1}, and
Theorem~\ref{T:2} follows immediately as well in light
of the facts that $Cg(5\cdot 10^{18})<7.06$ and $Cg(5\cdot 10^{55})<7$.

\vspace{2ex}
\noindent\textbf{Remark.}
It would be highly desirable to prove a form of Theorem~\ref{T:2} with a reasonable constant
when $p<10^{20}$.  For small $p$ the best result appears to be
due to A.~Brauer, using elementary methods.  In~\cite{bib:brauer.1932}, he shows that
$H\leq\sqrt{2p}+2$ for all $p$. 


\vspace{2ex}
\noindent\textbf{Acknowledgement.}
The author would like to thank Professor H. M. Stark for his
helpful suggestions.

\end{document}